\documentclass[11pt]{article}

\usepackage{amsmath,amssymb,amsthm,amscd,mathrsfs,eucal}
\usepackage{mathtools,stmaryrd}
\usepackage{latexsym}
\usepackage{xparse}
\usepackage{graphicx}
\usepackage[usenames]{color} 
\usepackage{enumerate}
\usepackage{mathscinet}
\usepackage{hyperref}
\usepackage{pgf,tikz}
\usetikzlibrary{arrows.meta}
\usepackage{hvlogos}
\usepackage{cleveref}
\usepackage[all,pdf]{xy}
\usepackage{lmodern}
\usepackage[affil-it]{authblk}
\usetikzlibrary{positioning}

\newtheorem{thm}[equation]{Theorem}

\newtheorem{lemma}[equation]{Lemma}
\newtheorem{prop}[equation]{Proposition}

\newtheorem*{conj*}{Conjecture}
\theoremstyle{definition}

\numberwithin{equation}{section}

\newcommand{\QQ}{\mathbb{Q}} 
\newcommand{\CC}{\mathbb{C}}  
\newcommand{\NN}{\mathbb{N}} 
\newcommand{\KK}{\mathbb{K}}

\newcommand{\qu}{\mathcal{O}}

\newcommand{\cZ}{\mathcal{Z}}

\newcommand{\Hf}{\mathcal{H}} 

\newcommand\fsl{\mathfrak{sl}}
\newcommand\End{\mathsf{End}}

\newcommand\rank{\mathsf{rank}}

\DeclareMathOperator{\GKdim}{GKdim}
\DeclareMathOperator{\kdim}{K.dim}
\DeclareMathOperator{\gfrac}{\mathrm{Frac}}
\DeclareMathOperator{\heit}{\mathrm{ht}}

\DeclareMathOperator{\Hom}{\mathrm{Hom}}
\DeclareMathOperator{\de}{\mathrm{d}}
\DeclareMathOperator{\si}{\mathrm{s}}

\DeclareMathOperator{\pois}{\mathsf{P}}
\DeclareMathOperator{\pdeg}{\mathrm{prim.deg}}
\DeclareMathOperator{\rdeg}{\mathrm{rat.deg}}
\DeclareMathOperator{\ldeg}{\mathrm{loc.deg}}

\DeclareMathOperator{\tr}{\mathrm{tr}}

\DeclareMathOperator{\n}{\mathfrak{n}}

\newcommand\g{\mathfrak{g}}
\newcommand\bor{\mathfrak{b}}

\newcommand\h{\mathfrak{h}}

\newcommand\ii{\mathfrak{i}}

\newcommand\W{\mathsf{W}}
\newcommand\w{\mathsf{w}}

\DeclareMathOperator{\spec}{\mathsf{Spec}}
\DeclareMathOperator{\hspec}{\Hf-\mathsf{Spec}}
\DeclareMathOperator{\Prim}{\mathsf{Prim}}

\usepackage[a4paper,margin=1.48in]{geometry}

\newcounter{marg}[section]

\usepackage[normalem]{ulem}
\newcommand\rsout{\bgroup\markoverwith{\textcolor{red}{\rule[0.5ex]{2pt}{0.4pt}}}\ULon}

\begin{document} 
\title{On the strong Dixmier--Moeglin equivalence for the quantum semisimple groups}
\author{Nirina Albert Razafimandimby\\ up202212304@edu.fc.up.pt}
\affil{Departamento de Matemática, Faculdade de Ciências, Universidade do Porto, \protect\\ Rua do Campo Alegre s/n, 4169--007 Porto, Portugal}
\date{}

	\maketitle
	
	\begin{description}
		\item[ Abstract:] 
		The Dixmier--Moeglin equivalence, when it holds for an algebra, characterizes its primitive ideals (that is,  the annihilators of irreducible representations) both algebraically and topologically within the prime spectrum. Bell, Launois, and Nolan \cite{BLN17} introduced a stronger version, known as the strong Dixmier--Moeglin equivalence. In this paper, we prove that the generic quantized coordinate ring $\qu_q(G)$ satisfies the strong Dixmier--Moeglin  equivalence for every connected, simply connected semisimple complex algebraic group $G$.
	\end{description}

\section{Introduction}
Throughout this paper, $\KK$ denotes an algebraically closed field of characteristic zero, every algebra is a unital associative $\KK$-algebra and every ideal is two-sided.

Classifying simple modules over an algebra is a fundamental and often difficult problem in representation theory. Simple modules play a role analogous to that of simple groups in group theory, as they are the building blocks of module categories. Their connection to ideals is captured by annihilators: the annihilator of a module consists of the elements of the algebra that act trivially on every element of that module. This naturally leads to the notion of primitive ideals, which are defined as the annihilators of simple left modules. While every primitive ideal is prime the converse is not always true: not all prime ideals are primitive.

The study of primitive ideals provides valuable insight into the structure and behavior of simple modules. In the 1980s, Dixmier and Moeglin introduced a characterization of primitive ideals in enveloping algebras of finite-dimensional Lie algebras in terms of algebraic and  topological properties of prime ideals. 

Let $A$ be a (possibly noncommutative) noetherian algebra, and let $\spec A$ denote its prime spectrum, which consists of all prime ideals of $A$,  equipped with the Jacobson--Zariski topology. The closed subsets in this topology are of the form  $V(I) = \{ P \in \spec A \mid P \supseteq I \}$ for an ideal $I$  of $A$. A prime ideal $P$ is said to be \textbf{locally closed} if the singleton $\{P\}$ is a locally closed subset;  that is, if it is the intersection of an open subset and a closed subset of $\spec A$. 
  Equivalently, $P$ is locally closed if and only if $P$ is strictly contained in the intersection of all prime ideals that strictly contain $P$.
  
Given a prime ideal $P \in \spec A$, we denote by $\gfrac(A/P)$ the Goldie quotient ring of the prime Noetherian ring $A/P$, which is obtained by localizing $A/P$ at the set of its regular elements. By Goldie's Theorem, $\gfrac(A/P)$ is a simple artinian ring. By the Artin-Wedderburn Theorem, the center of this ring is isomorphic to the center of the division ring $D$, i.e., $\cZ(\gfrac(A/P)) \simeq \cZ(D)$, which forms a field extension of the base field $\KK$. A prime ideal $P$ is called \textbf{rational} if this center, $\cZ(\gfrac(A/P))$, is an algebraic extension of $\KK$. In particular, as $\KK$ is algebraically closed, then every algebraic extension of $\KK$ is equal to $\KK$, and hence rationality is equivalent to $\cZ(\gfrac(A/P))= \KK$.
  Dixmier \cite{DI77} and Moeglin \cite{Mo80} proved that when $A$ is the enveloping algebra of a finite-dimensional complex Lie algebra, the notions of primitive, locally closed, and rational coincide. 
Irving and Small \cite{IS80} extended this result to an arbitrary base field of characteristic zero. This leads to the following definition. A noetherian algebra $A$ is said to satisfy the \textbf{Dixmier--Moeglin equivalence} if, for every prime ideal $P$ of $A$, the following conditions are equivalent: 
\begin{enumerate}
	\item[(1)] $P$ is locally closed;
	\item[(2)] $P$ is primitive;
	\item[(3)] $P$ is rational. 
\end{enumerate}
Recall that a noetherian algebra $A$ is said to satisfy the \textbf{Nullstellensatz} if  $\End_{A}(M)= \KK.1_{M}$ for every simple $A$-module $M$,  and if $A$ is a Jacobson ring  (meaning that every prime ideal is an intersection of primitive ideals). 

For a noetherian algebra satisfying the Nullstellensatz, one has the chain of implications:  $(1)\implies(2)\implies(3)$ (see \cite[Lemma II.7.15]{BG02}). In particular, these implications hold for countably generated noetherian algebras and when $\KK$ is uncountable (see \cite[II.7.16]{BG02}). Thus, the main difficulty in proving the Dixmier--Moeglin equivalence is establishing the implication: $(3)\implies(1)$. 

The Dixmier--Moeglin equivalence does not hold in all cases. Lorenz \cite{LO77} constructed an algebra with a primitive ideal that is not locally closed, while Irving \cite{IR79} provided an example of a rational ideal that is not primitive. Conversely, the equivalence has been established for various classes of algebras endowed with a Hopf structure, such as group algebras \cite{ZA71} and quantized coordinate rings \cite{GL00}. Motivated by these observations, Bell and Leung \cite[Conjecture 1.3]{BLE14} conjectured that every finitely generated complex Noetherian Hopf algebra of finite Gelfand--Kirillov dimension satisfies the Dixmier--Moeglin equivalence. As shown in \cite{BLSM15}, the equivalence can fail in the absence of such a structure; this suggests that the Hopf structure is critical to the conjecture, which remains open for the general case.

In 2017, Bell, Launois, and Nolan \cite{BLN17} introduced a stronger version of the Dixmier--Moeglin equivalence. They defined three degrees that measure how far a prime ideal is from being primitive, rational, or locally closed.

Let $A$ be a noetherian algebra and $P\in\spec A$. The \textbf{primitivity degree} of $P$ is
$$\pdeg P=\inf\{\heit Q\mid Q\in\Prim(A/P)\},$$
where $\Prim(A/P)$ is the set of primitive ideals of $A/P$. Thus, $P$ is primitive exactly when $\pdeg P=0$.

The \textbf{rationality degree} of $P$ is
$$\rdeg P=\tr.\deg_{\KK}\cZ\bigl(\gfrac(A/P)\bigr),$$
the transcendence degree over $\KK$ of the center of the quotient division ring of $A/P$. Thus, $P$ is rational exactly when $\rdeg P=0$.

The \textbf{local closure degree} of $P$ is the smallest nonnegative integer $d$ such that
$$\bigcap_{Q\in\spec_{>d}(A/P)}Q\neq\langle 0\rangle,$$
where $\spec_{>d}(A/P)$ consists of the prime ideals of $A/P$ having height greater than $d$. 
In particular, $P$ is locally closed exactly when $\ldeg P=0$.

The algebra $A$ is said to satisfy the \textbf{strong Dixmier--Moeglin equivalence} if
$$\pdeg P=\rdeg P=\ldeg P,
\quad\text{for every }P\in\spec A.$$
By contrast, $A$ satisfies the Dixmier--Moeglin equivalence precisely when
$$\ldeg P=0\Leftrightarrow\pdeg P=0\Leftrightarrow\rdeg P=0,
\quad\text{for every }P\in\spec A.$$
Thus, the strong Dixmier--Moeglin equivalence implies the Dixmier--Moeglin equivalence, but the converse need not hold. The enveloping algebra $U(\fsl_2)$ is a counterexample.

Examples of algebras known to satisfy the strong Dixmier--Moeglin equivalence include commutative affine domains, uniparameter quantum affine spaces and their corresponding tori, and quantum Schubert cells at non-roots of unity associated with finite-dimensional complex simple Lie algebras \cite{BLN17}. The equivalence has also been established for quantum special and general linear groups and for quantum matrices \cite[$\S$3.8]{NO17}.

The aim of this paper is to establish the strong Dixmier--Moeglin equivalence for generic quantized coordinate rings $\qu_q(G)$ with $G$ a connected, simply connected semisimple complex algebraic group. This class of algebras exhibits many properties that are analogous to those of classical coordinate rings $\qu(G)$, such as being noetherian \cite[Theorem I.8.18]{BG02} and satisfying the Nullstellensatz \cite[Corollary II.7.20]{BG02}. Our strategy is to use the Joseph parametrization of the prime spectrum of $\qu_q(G)$ by the double Weyl group $W\times W$ and the stratification theory of Goodearl and Letzter, along with combinatorial properties of the Weyl group. 

The paper is organized as follows. Section \ref{qsemisimpleg} introduces the basic definitions and constructions used throughout the paper and recalls Joseph's parametrization of $\spec \qu_q(G)$. Section \ref{sdmeqg} is devoted to our main results. In its first subsection, we establish the quasi-strong Dixmier--Moeglin equivalence, a weaker form of the strong Dixmier--Moeglin equivalence, as an intermediate step. Then, we present combinatorial results concerning Weyl groups. Finally, we prove the strong Dixmier--Moeglin equivalence for generic quantized coordinate rings $\qu_q(G)$, where $G$ is a connected, simply connected semisimple complex algebraic group.

\section{Quantum semisimple groups}\label{qsemisimpleg}
In this section, we will introduce the basic definitions and constructions that will be used throughout the rest of the paper. 

\subsection{Quantized enveloping algebras and their representations}
Let $G$ be a connected, simply connected complex semisimple algebraic group of rank $\ell$ and let $\g$ denote its Lie algebra. We write $\g= \n^{+}\oplus\h\oplus\n^{-}$ for a triangular decomposition of $\g$, where $\h$ is a Cartan subalgebra of $\g$, and $\n^{+}$ and $\n^{-}$ are the corresponding positive and negative nilpotent subalgebras, respectively. From this decomposition, we get that Borel subalgebras $\bor^{+}:= \h\oplus\n^{+}$ and $\bor^{-}:= \h\oplus\n^{-}$. Let $\pi=\{\alpha_1,\cdots,\alpha_\ell\}$ be the set of simple roots  and let $(-, -)$ be the bilinear form on $\h^\ast$ induced by the Killing form of $\g$. The (finite) Weyl group $W$ is generated by the simple reflections $s_1,\dots,s_\ell$ associated with the simple roots $\alpha_1,\dots,\alpha_\ell$.

Let $\pois(\pi)$ denote the weight lattice  of $\g$, and $\pois^{+}(\pi)$ denote the set of dominant integral weights. Recall that there is  a partial order denoted by $\le$ on $\pois(\pi)$, defined by  $\mu\le\lambda$ if $\lambda- \mu\in\pois^+(\pi)$.

For the remainder of the paper, let $q \in \KK^{\times}$ be a scalar that is not a root of unity. We denote by $U_q(\g)$ the quantized enveloping algebra of $\g$. Recall from \cite[2.1 and 3.1]{JL94} that this algebra is generated by $E_i, F_i,  K_i^{\pm1}$ with $1\leq i \leq \ell$. The reader is referred to \cite[2.1 and 3.1]{JL94} for the relations between these generators. As usual, we denote by $U_q(\bor^{+})$,  respectively, $U_q(\bor^{-})$, the subalgebra of $U_q(\g)$ generated by $E_1,\cdots, E_{\ell}, K_1^{\pm1},\cdots,K_{\ell}^{\pm1}$ respectively $F_1,\cdots,F_{\ell}, K_1^{\pm1},\cdots,K_{\ell}^{\pm1}$. 

Let $V$ be a $U_q(\g)$-module, and  let $\lambda\in\pois(\pi)$. Define
$$V_{\lambda}:= \{v\in V| K_i.v= q^{(\lambda,\alpha_i)}v,\quad
\forall 1\le i\le \ell\}.$$
If $V_{\lambda} \neq \{0\}$, we say that $\lambda$ is a weight of $V$.

Given a dominant integral weight $\lambda \in \pois^+(\pi)$, let $V(\lambda)$ denote the corresponding Verma module with highest weight vector $v_{\lambda}$. We refer to \cite[I.6.12]{BG02} for the construction and recall that the module is characterized by the following properties:
\begin{enumerate}
	\item[(i)]  $V(\lambda) = U_q(\g)v_{\lambda} = \bigoplus_{\mu \le \lambda} V(\lambda)_{\mu}$, where $V(\lambda)_{\mu}$ denotes the $\mu$-weight space of $V(\lambda)$;
	\item[(ii)] $V(\lambda)_{\lambda} = \KK v_{\lambda}$ and $E_i \cdot v_{\lambda} = 0$ for all $1 \le i \le \ell$. 
\end{enumerate}

\subsection{Quantized coordinate ring of $G$}

For $\lambda\in \pois^+(\pi), \xi\in V(\lambda)^{\ast}, v\in V(\lambda)$, define the matrix coefficients $c_{\xi,v}^{\lambda}\in(U_q(\g))^{\ast}$:
$$c_{\xi,v}^{\lambda}(x)= \xi(x.v), \quad \forall x\in U_q(\g).$$
The quantized coordinate ring of $G$ is the $\KK$-algebra generated by these matrix coefficients:
$$\qu_q(G):= \KK\langle c_{\xi,v}^{\lambda} | \lambda\in\pois^{+}(\pi), \xi\in V(\lambda)^{\ast}, v\in V(\lambda)\rangle.$$ 
It is well known, see for instance \cite[Proposition I.7.8]{BG02}, that 
$$\qu_q(G)= \KK\langle c_{\xi,v}^{\varpi_i} | i=1,\cdots,\ell , \xi\in V(\varpi_i)^{\ast},v\in V(\varpi_i)\rangle , $$
where $\varpi_1, \dots, \varpi_{\ell}$ are the fundamental weights.

  By  \cite[Lemma I.8.7]{BG02}, the quantized coordinate ring $\qu_q(G)$ admits a $\pois(\pi)\times\pois(\pi)$-grading:
\begin{equation}\label{decompqug}
	\qu_q(G)= \displaystyle\bigoplus_{(\beta,\eta)\in\pois(\pi)\times \pois(\pi)}\qu_q(G)_{\beta,\eta}, 
\end{equation} 
where
\begin{equation}\label{weigspaqug}
	\qu_q(G)_{\beta,\eta}:= \{y\in\qu_q(G)| y.K_i^{-1}= q^{(\beta,\alpha_i)}y, K_i.y= q^{(\eta,\alpha_i)}y, \forall 1\le i\le\ell\},
\end{equation}
for all $\beta,\eta\in\pois(\pi)$. From this decomposition, it follows that $\qu_q(G)$ is an integral domain; see \cite[Theorem I.8.9]{BG02}. Moreover,  $\qu_q(G)$ is a noetherian algebra by \cite[Theorem I.8.18]{BG02}. Thus $\qu_q(G)$ is a noetherian domain. 

\subsection{$\Hf$-stratification of the prime spectrum of $\qu_q(G)$}
Following \cite[I.7.1, II.1.18]{BG02}, let $\Hf$ be a maximal torus of $G$. We identify $\Hf$ with the algebraic torus $\Hom(\pois(\pi), \CC^{\times})$. The group $\Hf$ acts rationally on $\qu_q(G)$ via right winding automorphisms; specifically, for $\chi \in \Hf$ and $y \in \qu_q(G)_{\beta, \eta}$ (with $\beta, \eta \in \pois(\pi)$), the action is defined by $\chi \cdot y = \chi(\eta)y$. The rationality of this action is established in \cite[Exercise II.2.G]{BG02}.

An ideal $I \subseteq \qu_q(G)$ is called an $\Hf$-ideal if it is $\Hf$-stable, that is, $\chi(I) \subseteq I$ for all $\chi \in \Hf$. We say that an $\Hf$-ideal $P$ is $\Hf$-prime if, for any $\Hf$-ideals $I, J$ such that $IJ \subseteq P$, we have $I \subseteq P$ or $J \subseteq P$. Since $\qu_q(G)$ is noetherian, $\Hf$-prime ideals are precisely the prime $\Hf$-ideals; see \cite[Corollary II.1.12]{BG02}. We denote the set of $\Hf$-prime ideals of $\qu_q(G)$ by $\hspec\qu_q(G)$.  For a prime ideal $Q \in\spec\qu_q(G)$, define its $\Hf$-core by $$(Q : \Hf) = \bigcap_{h \in \Hf} h(Q).$$ 
This is the largest $\Hf$-stable ideal contained in $Q$, and it is an $\Hf$-prime ideal.

Consequently, one obtains the $\Hf$-stratification of the prime spectrum in the sense of Goodearl--Letzter \cite[2.1]{GL00}:
\begin{equation}\label{glhstrata}
	\spec\qu_q(G) = \bigsqcup_{J \in\hspec \qu_q(G)} \spec_J\qu_q(G),
\end{equation}
where  $$
\spec_J\qu_q(G)= \{Q\in\spec\qu_q(G): (Q: \Hf)= J\}.$$
Thus, the stratum associated with $J$ consists precisely of those prime ideals whose $\Hf$-core is  $J$.


	\subsection{Joseph's stratification of the prime spectrum of $\qu_q(G)$}\label{specqg}
	
	For $\mu\in\pois^{+}(\pi)$, let $\Omega(V(\mu))$ denote the collection of weights of $V(\mu)$,  counted with their multiplicities. Given weights $\lambda,\eta\in\Omega(V(\mu))$ and vectors $\xi\in V(\mu)^{\ast}_{\lambda}, v\in V(\mu)_{\eta}$, we will, by abuse of notation, write $c_{\lambda,\eta}^{\mu}:= c_{\xi,v}^{\mu}$ . 
	
The Weyl group $W$ is a Coxeter group, so each element $w \in W$ has a well-defined length $\ell(w)$. Specifically, $\ell(w)$ is the smallest non-negative integer $r$ such that $w$ can be expressed as a product of $r$ simple reflections: $w = s_{i_1}s_{i_2}\cdots s_{i_r}$. Since $W$ is finite, it has a unique element of maximal length, denoted by $w_0$. Set $\W= W\times W$. We write its elements as $\w= (w_{+},w_{-})$ with $w_+, w_-\in W$.
	
For $\lambda \in\pois(\pi)$ and $w \in W$, it follows from \cite[4.3.6(iv)]{JO95} that $\dim_{\KK}(V_{w\lambda}(\lambda)) = \dim_{\KK}(V(\lambda)^*_{w\lambda}) = 1$. We thus select nonzero vectors $v_{w\lambda} \in V(\lambda)_{w\lambda}$ and $\xi_{w\lambda} \in V(\lambda)^*_{w\lambda}$ such that $\xi_{w\lambda}(v_{w\lambda}) = 1$.
	
	Let $M$ be a $U_q(\g)$-module and let $N\subseteq M$. The annihilator, or orthogonal complement of $N$ in $M^{\ast}$ is:
	$$N^{\perp}:= \{f\in M^{\ast}| N\subseteq\ker(f)\}.$$
	For $w\in W$ and $\lambda\in\pois^{+}(\pi)$, define
	\begin{align*}
		V_{w}^{+}(\lambda)^{\perp}:&= \{c_{\xi,v}^{\lambda}| v\in V(\lambda)_{\lambda}, \xi\in(U_q(\bor^{+})v_{w\lambda})^{\perp}\subseteq V(\lambda)^{\ast}\},\\
		V_{w}^{-}(\lambda)^{\perp}:&= \{c_{\xi,v}^{\lambda}| v\in V(\lambda)_{w_0\lambda}, \xi\in(U_q(\bor^{-})v_{ww_0\lambda})^{\perp}\subseteq V(\lambda)^{\ast}\}.
	\end{align*}
	We then define $J_{w}^{\pm}:= \sum_{\lambda\in\pois^{+}(\pi)}\qu_q(G)V_{w}^{\pm}(\lambda)^{\perp}$. For $\w= (w_{+},w_{-})\in\W$, set
	\begin{equation}\label{minideal}
		J_{\w}:= J_{w_{+}}^{+}+ J_{w_{-}}^{-}.
	\end{equation}
	
	By \cite[Exercise II.4.C]{BG02}, the matrix coefficients $c_{\xi,v}^{\lambda}$ are $\Hf$-eigenvectors. Consequently, each $J_w^{\pm}$, and hence $J_{\w}$ is  $\Hf$-stable under the $\Hf$-action. Therefore, $J_{\w}$ is an $\Hf$-ideal of $\qu_q(G)$. 
	
	For $\mu\in\pois^{+}(\pi)$ and  $w\in W$, define
	$$c_{w\mu}:= c_{w\mu, \mu}^{\mu}, \quad d_{w\mu}:= c_{ww_0\mu,w_0\mu}^{\mu}.$$
	
    By \cite[9.1.5]{JO95}, the matrix coefficients  $c_{w\lambda}$ and $d_{w\mu}$ are normal modulo $J_{w}^+$ and  $J_w^-$, respectively.
	
	For $\lambda\in\pois^+(\pi)$ and $\w=(w_+,w_-)\in \W$, set $\lambda_0:= -w_0\lambda$, and define
	$$b_{\w,\lambda}:= c_{w_+\lambda}d_{w_-\lambda_0}.$$ 
	
	Let $E_{\w}$ be the multiplicative subset of $\qu_q(G)/J_{\w}$ generated by the images of the elements $b_{\w,\lambda}$, for $\lambda\in\pois^+(\pi)$ together with $\KK^{\times}$. Since $J_{w_{\pm}}^{\pm}\subseteq J_{\w}$, their images are normal modulo $J_{\w}$. Furthermore,  by \cite[9.1.5]{JO95} the factors $c_{w_+\lambda}$ and $d_{w_-\lambda_0}$ commute up to nonzero scalars  modulo $J_{\w}$.  Hence $b_{\w,\lambda}$ is normal in $\qu_q(G)/J_{\w}$.  Using (\cite[9.1.10]{JO95}, we also obtain: $b_{\w,\lambda}b_{\w,\mu}\equiv b_{\w,\lambda+\mu}\mod J_{\w}$, up to a nonzero scalar multiple, for all  $\lambda,\mu\in\pois^{+}(\pi)$. Consequently, the multiplicative system $E_{\w}$ is generated, up to nonzero scalar multiples, by the images in $\qu_q(G)/J_{\w}$ of the finite set $\{ b_{\w,\alpha}, \alpha\in\pi\}$.  Since these generators are normal $\Hf$-eigenvectors and commute pairwise up to nonzero scalar multiples, $E_{\w}$ is an Ore set. We may therefore form the localization
	$$R_{\w}:= (\qu_q(G)/J_{\w})E_{\w}^{-1}.$$
	
	As noted above $J_{\w}$ is an $\Hf$-ideal, and $E_{\w}$ is generated by $\Hf$-eigenvectors. Therefore, the $\Hf$-action on $\qu_q(G)$ extends to a rational action on $R_{\w}$ by $\KK$-algebra automorphisms. Moreover, for every $\w\in\W$, the algebra $R_{\w}$ is nonzero and $\Hf$-simple, that is, its only $\Hf$-stable ideals  are $\langle 0 \rangle$ and $R_{\w}$. These fact follow from \cite[10.3.2]{JO95} and \cite[Corollary II.4.10]{BG02}. 
	
	For each $\w\in\W$, define 
	$$\spec_{\w}\qu_q(G):= \{P\in\spec\qu_q(G)| P\supseteq J_{\w} ~~ \text{and $(P/J_{\w})\cap E_{\w}= \emptyset$}\} $$
	and 
	$$\Prim_{\w}\qu_q(G)= \Prim\qu_q(G)\cap\spec_{\w}\qu_q(G).$$ 
These sets partition $\spec\qu_q(G)$ and $\Prim\qu_q(G)$, respectively, as the following result shows.
	\begin{prop}\label{hspecqug}
		\begin{itemize}
			\item[(a)] The set of $\Hf$-prime ideals of $\qu_q(G)$ is $$\hspec\qu_q(G)= \{J_{\w}| \w\in\W\}.$$
			\item[(b)] For every $w$ in $\W$, $\spec_{J_{\w}}(\qu_q(G))$ is equal to $\spec_{\w}(\qu_q(G))$. Consequently, the partition of $\spec(\qu_q(G))$ into these subsets coincides with the $\Hf$-stratification in (\ref{glhstrata}). 
			\item[(c)] The set of primitive ideals in each stratum $\spec_{\w}\qu_q(G)$ coincides with the set of maximal elements of $\spec_{\w} \qu_q(G)$. 
		\end{itemize}
	\end{prop}
	\begin{proof}
		Parts (a) and (b) follow from \cite[Theorem 2.1] {YA14}; see also \cite[Proposition II.4.11, Corollary II.4.12]{BG02}. Part (c) follows from part (b) together with the standard result on $\Hf$-stratifications stating that the primitive ideals are exactly the prime ideals that are maximal, with respect to inclusion, within their $\Hf$-strata; see \cite[Corollary 6.9]{GL00}. 
	\end{proof}
	
	Note for later use that, since $E_{\w}$ is an Ore set, extension ($P\mapsto \frac{P}{J_{\w}}R_{\w}$) and contraction ($\frac{Q}{J_{\w}}\mapsto Q\cap\qu_q(G)$) induce mutually inverse homeomorphisms between $\spec_{\w}\qu_q(G)$ and $\spec R_{\w}$.

	\section{Strong Dixmier--Moeglin equivalence for $\qu_q(G)$}\label{sdmeqg}
	In this section, we prove the strong Dixmier--Moeglin equivalence for the quantized coordinate ring $\qu_q(G)$. We first establish a weaker version of the strong Dixmier–Moeglin equivalence, namely the \textbf{quasi-strong Dixmier--Moeglin equivalence}. More precisely, we prove that
	$$\ldeg P= \rdeg P,\quad\forall P\in \spec\qu_q(G).$$
	Then, we prove some necessary combinatorial properties of Weyl groups before showing that the primitivity and the local closure degrees coincide, thus establishing our main result. 	
	\subsection{Quasi-strong Dixmier--Moeglin equivalence for $\qu_q(G)$}\label{qsdmeqg}	
	
	\begin{prop}\label{qsdmequg}
		Let $G$ be a connected, simply connected complex semisimple algebraic group. Then, the quantized coordinate ring $\qu_q(G)$ satisfies the quasi-strong Dixmier--Moeglin equivalence. 
	\end{prop}
	\begin{proof}
	For $\lambda\in \pois^+(\pi)$, let $c_{\lambda,\lambda}^{\lambda}:= c_{\xi_{\lambda},v_{\lambda}}^{\lambda}$. It follows from \cite[ 9.1.10 and 9.3.10]{JO95} that $\Sigma:= \{c_{\lambda,\lambda}^{\lambda}| \lambda\in\pois^+(\pi)\}$ is an Ore set of regular elements in $\qu_q(G)$. Moreover, for every $P\in\spec\qu_q(G)$, we have $\Sigma\cap P= \emptyset$. 
	
	It follows from  \cite[Lemmas 6.4.1 and 6.4.2]{Ca03} that $$\qu_q(G)\Sigma^{-1}\simeq BS^{-1},$$
	where $B$ is a uniparameter quantum nilpotent algebra in the sense of
	\cite[Definition 3.1]{LLR06}, and $S$ is an Ore set generated by regular
	normal elements.
	
		By \cite[Theorem 7.1]{BLN17}, $B$ satisfies the quasi-strong Dixmier--Moeglin equivalence, and then, by \cite[Lemma 2.3]{BLN17}, so does the localization $BS^{-1}$. This shows that  the algebra $\qu_q(G)\Sigma^{-1}$ satisfies the quasi-strong Dixmier--Moeglin equivalence. By construction, $\Sigma$ is disjoint from every prime ideal of $\qu_q(G)$. Consequently, $\spec\qu_q(G)= \{P\in\spec\qu_q(G)| P\cap \Sigma= \emptyset\}$. It follows by \cite[Theorem 10.20]{GW04} that  extension $(P\mapsto P\Sigma^{-1})$ and contraction $(Q\mapsto Q\cap\qu_q(G))$ induce mutually inverse, order-preserving homeomorphisms between $\spec\qu_q(G)$ and $\spec\qu_q(G)\Sigma^{-1}$. Let $P$ be a prime ideal of $\qu_q(G)$. Then $P\Sigma^{-1}$ is a prime ideal of $\qu_q(G)\Sigma^{-1}$ and there is a natural isomorphism
		\begin{equation}\label{quoqugiso}
			\qu_q(G)\Sigma^{-1}/P\Sigma^{-1}\simeq(\qu_q(G)/P)\overline{\Sigma}^{-1}, 
		\end{equation}
		where $\overline{\Sigma}$ denotes the image of $\Sigma$ in $\qu_q(G)/P$. Moreover, $\spec\qu_q(G)/P$ is homeomorphic to the closed subset $\{Q\in\spec\qu_q(G)| Q\supseteq P\}$ of $\spec\qu_q(G)$. Hence, localization induces a homeomorphism $\spec(\qu_q(G)/P)\overline{\Sigma}^{-1}\simeq\spec\qu_q(G)/P$. Therefore,
		\begin{align*}
			\ldeg P= \ldeg\langle 0 \rangle_{\qu_q(G)/P}= \ldeg\langle 0\rangle_{(\qu_q(G)/P)\overline{\Sigma}^{-1}},
		\end{align*}
		where $\langle 0 \rangle_{\qu_q(G)/P}$ and $\langle 0 \rangle_{(\qu_q(G)/P)\overline{\Sigma}^{-1}}$ denote the zero ideals of $\qu_q(G)/P$ and $(\qu_q(G)/P)\overline{\Sigma}^{-1}$, respectively. It follows from the isomorphism (\ref{quoqugiso}) that \begin{align*}
			\ldeg P&= \ldeg\langle 0 \rangle_{(\qu_q(G)/P)\overline{\Sigma}^{-1}}= \ldeg P\Sigma^{-1}\\
			&= \rdeg P\Sigma^{-1}\\
			&= \rdeg P.
		\end{align*}
		The second to last equality results from the fact that $\qu_q(G) \Sigma^{-1}$ satisfies the quasi-strong Dixmier--Moeglin equivalence, while the last equality follows from the isomorphism (\ref{quoqugiso}). Therefore, $\qu_q(G)$ satisfies the quasi-strong Dixmier--Moeglin equivalence.
	\end{proof}

	\subsection{Combinatorial properties of Weyl groups}\label{combWeylgroup}
	To establish the strong Dixmier--Moeglin equivalence, we will use
	several combinatorial properties of the Weyl group $W$. We begin by	defining
	$$T:= \{ws_{\alpha}w^{-1}|  w\in W, \alpha\in\pi\},$$
	the set of reflections in $W$. 
	
	Following \cite[Definition 2.1.1 ]{BB05}, for $u,w\in W$, we write $u\to w$ if $u^{-1}w\in T$ and $\ell(u)< \ell(w)$. The Bruhat order  on $W$, is the partial order generated by this relation. Thus, write $u\le w$ if either $u= w$ or there exists a sequence
	$$u= u_0\to u_1\to \cdots\to u_{k-1}\to u_k= w.$$ 
	of elements of $W$.

	For $w\in W$, let $\si(w)$ denote the minimum number of reflections required to express $w$ as their product, where the reflections need not be simple. We have the following standard properties. 
	\begin{lemma}\label{metriconw} For all $u, w\in W$:
		\begin{enumerate}
			\item[(a)] $\si(w)= \si(w^{-1})$,
			\item[(b)] The function $\si(u^{-1}w)$ defines a metric on $W$. 
		\end{enumerate}
	\end{lemma}
	\begin{proof}
		\begin{enumerate}
			\item[(a)] Let $k= \si(w)$. Then $w$ can be written as a product of $k$ reflections: $w= t_1t_2\cdots t_k$. Taking inverses gives: $w^{-1}= t_k\cdots t_2t_1$, since every reflection is an involution. Hence $\si(w^{-1})\le\si(w)$. Applying the same argument to $w^{-1}$ yields $\si(w)\le\si(w^{-1})$. Therefore, $\si(w)= \si(w^{-1})$. 
			\item[(b)] Define $\de(u, w):= \si(u^{-1}w)$. We verify the metric axioms. Since only the identity element can be expressed as a product of zero reflections, we have $\de(u, w)= 0\iff \si(u^{-1}w)= 0 \iff u^{-1}w= e\iff u= w$. Moreover, by part (a), $\de(u,w)= \de(w, u)$.
			
			It remains to verify the triangle inequality. Let $u, v, w\in W$, and choose minimal expressions $u^{-1}v= t_1t_2\cdots t_k$ and $v^{-1}w= t_1't_2'\cdots t_m'$, where each $t_i$ and $t'_{j}$ is a reflection. Then, $$u^{-1}w= u^{-1}v(v^{-1}w)= t_1t_2\cdots t_kt_1't_2'\cdots t_m'.$$ 
			Thus $u^{-1}w$ is a product of $k+ m$ reflections, and hence $\si(u^{-1}w)\leq k+ m= \si(u^{-1}v)+ \si(v^{-1}w)$. Therefore, 
			$$\de(u, w)\leq \de(u,v)+ \de(v, w).$$
			Hence $\de(u,w)= \si(u^{-1}w)$ defines a metric on $W$. 
		\end{enumerate}
	\end{proof}
	
	\begin{lemma}\label{smalllenghtpro} Let $u, u_1, u_2, w_1, w_2\in W$. Then 
		\begin{enumerate}
			\item[(a)] If $u\le w$, then $\ell(u)+ \si(u^{-1}w)\le\ell(w)$.
			\item[(b)] If $u_1\le w_1$ and $u_2\le w_2$, then 
			$$\ell(u_1)+ \ell(u_2)+ \si(u_1^{-1}u_2)\le\ell(w_1)+ \ell(w_2)+ \si(w_1^{-1}w_2).$$
		\end{enumerate} 
	\end{lemma}
	\begin{proof}
		\begin{enumerate}
			\item[(a)] 
				We proceed by induction on $k = \si(u^{-1}w)$. 
				The base case $k=0$ implies $u=w$, which gives $\ell(u) + 0 = \ell(w)$, satisfying the inequality.
				
				Assume $k \ge 1$ and that the result holds for all elements with reflection length strictly less than $k$. 
				By the properties of the Bruhat order \cite[Theorem 2.2.6]{BB05}, since $u< w$, there exists a saturated chain $u = u_0 < u_1 < \dots < u_k = w$ such that $\ell(u_i) = \ell(u) + i$ for all $i \in \{1, \dots, k\}$. 
				For each step $i \in \{1, \dots, k\}$, the property $\ell(u_{i-1}) < \ell(u_i)$ implies there exists a reflection $t_i \in T$ such that $u_i = u_{i-1}t_i$ and $\ell(u_i) = \ell(u_{i-1}) + 1$. 
				By composing these relations, we obtain the factorization $w = u t_1 t_2 \dots t_k$, which implies $u^{-1}w = t_1 t_2 \dots t_k$.
				
				Define $u' = u_1 = u t_1$. Since $u < u_1$, we have $\ell(u') = \ell(u) + 1$. 
				Consider the element $(u')^{-1}w = t_2 \dots t_k$, which is a product of $k-1$ reflections. This implies $\si((u')^{-1}w) \le k - 1$. 
				Since $u' \le w$ in the Bruhat order, we apply the inductive hypothesis to $u'$ and $w$:
				 $$\ell(w) \ge \ell(u') + \si((u')^{-1}w)$$
				Substituting $\ell(u') = \ell(u) + 1$ and $\si((u')^{-1}w) = k-1$, we obtain:
				\[ \ell(w) \ge (\ell(u) + 1) + (k - 1) = \ell(u) + k \]
				Thus, $\ell(w) \ge \ell(u) + \si(u^{-1}w)$, completing the inductive step.

			\item [(b)] We use the triangular inequality from Lemma \ref{metriconw}(b):
			$$\si(u_1^{-1}u_2)\le\si(u_1^{-1}w_1)+ \si(w_1^{-1}w_2)+ \si(w_2^{-1}u_2).$$
			This can be rewritten using the equality from Lemma \ref{metriconw}(a):
			$$\si(u_1^{-1}u_2)\le\si(u_1^{-1}w_1)+ \si(w_1^{-1}w_2)+ \si(u_2^{-1}w_2).$$
			Then we can add lengths:
			\begin{align*}
				\ell(u_1)+ \ell(u_2)+ \si(u_1^{-1}u_2)\le\ell(u_1)+ \si(u_1^{-1}w_1)+ \ell(u_2)+ \si(u_2^{-1}w_2)+ \si(w_1^{-1}w_2).
			\end{align*}
			By applying (a) to the relations  $u_1\leq w_1$ and $u_2\leq w_2$:
			$$	\ell(u_1)+ \ell(u_2)+ \si(u_1^{-1}u_2)\le \ell(w_1)+ \ell(w_2)+ \si(w_1^{-1}w_2).$$
			Hence the result holds.
		\end{enumerate}
	\end{proof}
	
	\subsection{Strong Dixmier--Moeglin equivalence for $\qu_q(G)$}\label{spdmeOqG}
 Throughout this section, let $G$ be a connected, simply connected complex semisimple algebraic group, and let $\qu_q(G)$ be its quantized coordinate ring, where $q \in \KK^{\times}$ is transcendental over $\QQ$. Although we previously assumed only that $q$ is not a root of unity, we now assume that $q$ is transcendental over $\QQ$ to align with the literature regarding the poset structure. Note that while Lemma \ref{hspecqug}(a) establishes a bijection between $\hspec\qu_q(G)$ and $\W$ under the weaker assumption that $q$ is not a root of unity, the order-reversing poset isomorphism between $\hspec\qu_q(G)$ and $\W$ used in the proof of Proposition \ref{pairineqh} is provided by Joseph \cite{JO95}, which requires the stronger assumption that $q$ is transcendental.  For $\mu\in\pois(\pi)$ and $\w= (w_+,w_-)\in\W$, define $a_{\w,\mu}:= c_{w_+\mu}^{-1}d_{w_-\mu_0}\in R_{\w}$. These elements satisfy $a_{\w,\mu}a_{\w,\nu}= a_{\w,\mu+ \nu}$ up to non-zero scalars. Let $\cZ_{\w}$ denote the center of $R_{\w}$.

 \begin{prop}\label{centstrat}[\cite[Theorem II.4.14]{BG02}, \cite[Theorems 3.1 and 2.3(ii)]{YA14}]
 	Let $\w= (w_+,w_-)\in\W$, and define 
 	$$\pois^{\w}(\pi):= \{\mu\in\pois(\pi): w_+\mu- w_-\mu=0\}.$$
 	 This is a free abelian subgroup of $\pois(\pi)$. Let $\mu_1,\cdots,\mu_r$ be a basis for  $\pois^{\w}(\pi)$. Then the following hold:
 	   \begin{enumerate}
 		\item[(a)] The center $\cZ_{\w}$ is a Laurent polynomial ring:
 		$$\cZ_{\w}= \KK[a_{\w,\mu_1}^{\pm1},\cdots,a_{\w,\mu_r}^{\pm1}]$$
 		where $r= \rank(G)- \si(w_-^{-1}w_+)$. 
 		\item[(b)] The stratum $\spec_{\w}\qu_q(G)$  is homeomorphic to $\spec\cZ_{\w}$ via the map $\ii_{\w}: \spec\cZ_{\w}\to\spec_{\w}\qu_q(G)$ defined as follows. For each $Q_0\in\spec \cZ_{\w}$, the image $\ii_{\w}(Q_0)$ is the unique ideal of $\qu_q(G)$ containing $J_{\w}$ such that
 		$$\ii_{\w}(Q_0)/J_{\w}= (Q_0R_{\w})\cap(\qu_q(G)/J_{\w}).$$
 	\end{enumerate}
 \end{prop}

Motivated by Proposition \ref{centstrat}(b), we define the Krull dimension of each stratum $\spec_{\w}\qu_q(G)$ to be the Krull dimension of the corresponding center $\cZ_{\w}$. Throughout the remainder of this section, we equip $\W= W\times W$ with the product order, so  $(w_{+},w_-) \le (w_+',w_{-}')$  if and only if $ w_+ \le w_+'$ and $w_- \le w_-'$.

  	\begin{prop}\label{pairineqh}
  Let	$J_{\w}\subseteq J_{\w'}$ be a pair of $\Hf$-prime ideals  of $\qu_q(G)$. Then
  		$$\kdim(\spec_{\w}\qu_q(G))+ \heit J_{\w}\le\kdim(\spec_{\w'}\qu_q(G))+ \heit J_{\w'}.$$
  	\end{prop}
  	
  	\begin{proof}
  		For each $\w= (w_+,w_-)\in\W$, define $\ell(\w):=\ell(w_+)+ \ell(w_-)$ and $\si(\w):= \si(w_-^{-1}w_+)$.  For any $P\in\spec_{\w}(\qu_q(G))$,  \cite[Proposition 10.3.9]{JO95} gives
  		\begin{equation}\label{gdimquo}
  			\GKdim(\qu_q(G)/P)= \ell(\w)+ \si(\w)+ \kdim(\cZ_{\w}/(\cZ_{\w}\cap(P/J_{\w})E_{\w}^{-1})).
  		\end{equation}
  		In particular, taking $P= J_{\w}$ and applying Proposition \ref{centstrat}(a), we obtain
  		\begin{equation}
  			\GKdim(\qu_q(G)/J_{\w})= \ell(\w)+ \rank(G). 
  		\end{equation}
  		Tauvel's height formula \cite[Theorem 9.4]{YA14} yields
  		\begin{equation}\label{tauvelheit}
  			\GKdim(\qu_q(G)/J_{\w})+ \heit J_{\w}= \GKdim(\qu_q(G)). 
  		\end{equation}
  		
  		The same theorem gives $\GKdim(\qu_q(G))= \dim(G).$
        Hence
  		 \begin{equation}\label{heightjw}
  		 	\heit J_{\w}= \dim(G)- \rank(G)- \ell(\w). 
  		 \end{equation}
  		 On the other hand, Proposition \ref{centstrat}(a) implies that 
  		 $$\kdim(\spec_{\w}\qu_q(G))= \rank(G)- \si(\w).$$ 
  		 Combining this equality with (\ref{heightjw}), we obtain
  		 \begin{equation}\label{dimstratheit}
  		 	\kdim(\spec_{\w}\qu_q(G))+ \heit J_{\w}= \dim(G)- \ell(\w)- \si(\w). 
  		 \end{equation}
  		 Now suppose that $J_{\w}\subseteq J_{\w'}$. We write $\w= (w_+,w_-)$ and $\w'= (w_+',w_-')$. In particular, \(J_{w_+}^{+}\subseteq J_{w'_+}^{+}\). Hence, by
  		 \cite[Proposition 4.4.5 and the proof of Lemma 10.3.6]{JO95}, we have $V_{w_+}(\lambda)^{\perp}\subseteq V_{w'_+}(\lambda)^{\perp}$ for every $\lambda\in\pois^{+}(\pi)$, which implies that $w'_+\leq w_+$. Similarly, $w'_-\leq w_-$.
  		  Therefore, Lemma \ref{smalllenghtpro}(b) gives
  		\begin{equation}\label{lengthinq}
  			\ell(\w')+ \si(\w')\le\ell(\w)+ \si(\w).
  		\end{equation}
   Equivalently,
  		\begin{equation}\label{dimstratineq}
  			\dim(G)- \ell(\w)- \si(\w)\le\dim(G)- \ell(\w')- \si(\w'). 
  		\end{equation}
  		Using (\ref{dimstratheit})  for $\w$ and $\w'$, we conclude that
  	\[
  	\kdim(\spec_{\w}\qu_q(G)) + \heit J_{\w} \le \kdim(\spec_{\w'}\qu_q(G)) + \heit J_{\w'}.
  	\]
  	\end{proof}
We are now ready to establish the main result.
  	\begin{thm}
  		Let $G$ be a  connected, simply connected complex semisimple algebraic group $G$ and let  $q\in\KK^{\times}$ be transcendental over $\QQ$. Then the quantized coordinate ring $\qu_q(G)$ satisfies the strong Dixmier--Moeglin equivalence. 
  	\end{thm}
  	
  	\begin{proof}
  		To prove the theorem, we use the following foundational results:
  		\begin{enumerate}
  			\item The algebraic torus $\Hf$ acts rationally on $\qu_q(G)$ by $\KK$-algebra automorphisms; see \cite[Exercise II.2.G]{BG02}.
  			\item  The set $\hspec\qu_q(G)$ is finite; see Proposition \ref{hspecqug}(a).
  			\item By \cite[Proposition II.8.3]{BG02}, every $J\in\hspec\qu_q(G)$ is strongly $\Hf$-rational. That is, the subfield of $\Hf$-fixed elements of  $\cZ(\gfrac(\qu_q(G)/J))$ is $\KK$. 
  		\end{enumerate}
  		
Furthermore, $\qu_q(G)$ satisfies the Nullstellensatz over $\KK$ by \cite[Corollary II.7.20]{BG02}, and its prime spectrum is catenary. Therefore, Proposition \ref{pairineqh} shows that  $\qu_q(G)$ satisfies all the hypotheses of \cite[Proposition 5.2]{BLN17}. Consequently, for any $P\in\spec_{\w}\qu_q(G)$ and any primitive ideal  $M\supseteq P$ belonging to $\spec_{\w}\qu_q(G)$, we have
  		\begin{equation}\label{pdegqug}
  			\pdeg P=\heit M/P=  \kdim(\spec_{\w}\qu_q(G))+ \heit J_{\w}- \heit P.
  		\end{equation}
  		We next claim that if $Q\in\spec\qu_q(G)\setminus\spec_{\w}\qu_q(G)$ such that $Q\supseteq P$, then $b_{\w}\in Q$, where $b_{\w}= b_{\w,\alpha_1}\cdots b_{\w,\alpha_{\ell}}$. Indeed, since $Q\notin\spec_{\w}\qu_q(G)$, then  $(Q/J_{\w})\cap E_{\w}\neq\emptyset$. Every element of $E_{\w}$ has the form $kb_{\w,\alpha_1}^{n_1}b_{\w,\alpha_2}^{n_2}\cdots b_{\w,\alpha_{\ell}}^{n_{\ell}}+ J_{\w}$, where $k\in\KK^{\times}$ and $n_1,\cdots,n_{\ell}\in\NN$.  Thus, for some such $k, n_1,\cdots, n_{\ell}$, we have $kb_{\w,\alpha_1}^{n_1}b_{\w,\alpha_2}^{n_2}\cdots b_{\w,\alpha_{\ell}}^{n_{\ell}}\in Q$. By \cite[Théorème 6.4.1(1)]{Ca03}, $Q$ is completely prime. Hence, $b_{\w,\alpha_i}\in Q$ for some $i$, and therefore $b_{\w}\in Q$. This proves the claim.

  		  Now, let $Q\in\spec\qu_q(G)$ satisfy $Q\supseteq P$ and $\heit (Q/P)> \pdeg P$. It follows from (\ref{pdegqug}) that  $Q\notin\spec_{\w}\qu_q(G)$. By the claim, $b_{\w}\in Q$. On the other hand, $(P/J_{\w})\cap E_{\w}= \emptyset$, so $b_{\w}\notin P$. Consequently, 
  		  $$P\subsetneq \bigcap_{\substack{P\subseteq Q\in\spec\qu_q(G)\\\heit(Q/P)> \pdeg P}}Q.$$
  		  It follows that $\ldeg P\leq\pdeg P$. The reverse inequality is know by \cite[Proposition 2.1]{BLN17}. Hence
  		   \begin{equation}\label{ldegispdeg}
  		  	\ldeg P=  \pdeg P.
  		  \end{equation}
  		   Combining (\ref{ldegispdeg}) with Proposition \ref{qsdmequg} completes the proof.
  	\end{proof}

\section*{Acknowledgments}

We acknowledge the support from the bilateral cooperation French-Portuguese PESSOA Program (refs. PHC PESSOA 2025 - 54952PC and 2025.08159.CBM). This work was also partially supported by CMUP (under project \href{https://doi.org/10.54499/UID/00144/2025}{UID/00144/2025}) and an FCT scholarship (ref. UI/BD/154942/2023).

The author is fully responsible for the mathematical content, logical arguments, and citations. AI tools were only used for English language refinement.

\bibliographystyle{abbrv}
\bibliography{biblio.bib}
\end{document}